\documentstyle[12pt]{article}
\let\text\mbox

\newcommand\qed{\hfill $\diamond$ \vspace{5mm}}
\newcommand\ptwo{\text{\bf P\/}^2}
\newcommand\pthree{\text{\bf P\/}^3}
\newcommand\pn{\text{\bf P\/}^n}
\newcommand\C{\text{\bf C}}

\newcommand\gi{\text{gin}(I)}

\newcommand\giG{\text{gin}(I_{\Gamma})}
\newcommand\G{\Gamma}
\newcommand\lm{\lambda}
\newcommand\al{\alpha}
\newcommand\be{\beta}

\newcommand\ga{\gamma}
\newcommand\p{\partial}

\begin{document}

\title{The Connectedness of Space Curve Invariants}
\author{Michele Cook}
\maketitle

\newtheorem{theorem}{Theorem}
\newtheorem{lemma}[theorem]{Lemma}
\newtheorem{corollary}[theorem]{Corollary}
\newtheorem{proposition}[theorem]{Proposition}
\hyphenation{fam-ily homo-gen-eous poly-nom-ials gen-er-ated}

\section*{Introduction}
Let $ S = \C[x_1, \dots , x_{n+1}] $ be
the set of polynomials in n+1 variables over {\bf C} (usually corresponding
to the homogeneous coordinate ring of $\pn$).
Let $\succ$ be the reverse lexicographical order on the monomials of $S$.
Given a homogeneous ideal $I$ in $S$, $\gi$, the
{\bf generic initial ideal} of $I$, is the Borel-fixed
monomial ideal associated to $I$.
(For details see [B], [BM], [BS] or [Gr].)

The generic initial ideal reflects much of the structure of the
original ideal, for example, it has the same Hilbert function
and the same regularity.
Also, many geometric problems can be reduced to combinatorial
ones which may be phrased in terms of generic initial ideals.
Thus, we would like to have some rules governing the kinds of
Borel-fixed monomial ideals which can occur from geometry.

The first set of rules are supplied by a result of
Gruson and Peskine ([GP]):

\begin{theorem}[Gruson-Peskine]
Let $\G \subset \ptwo$ be a set of points in general position,
with generic initial ideal
$\giG = (x_1^s, \ x_1^{s-1}x_2^{\lm_{s-1}}, \ \dots , \
x_1x_2^{\lm_{1}},\ x_2^{\lm_{0}}).$
Then the {\bf invariants}, $\{\lm_{i}\}_{i=0}^{s-1}$, satisfy
$$ \lm_{i+1} + 2 \geq \lm_i \geq \lm_{i+1} + 1 \ \
\text{for all} \ i <  s-1  $$
and we say that the invariants of $\G$ are {\bf connected}.
\end{theorem}

Associated to a space curve there are invariants which generalize the
invariants of points in $\ptwo$ (see below).
The main result of this paper
is to show that the invariants of a
reduced, irreducible, non-degenerate curves in $\pthree$ also
satisfy a connectedness property:

Let $C$ be a reduced, irreducible, non-degenerate curve in $\pthree$,
As $I_C$ is saturated, the generators of $\gi$ will be of the form
$x_1^{i_1}x_2^{i_2}x_3^{f(i_1,i_2)}$.
\newline
Let
$$  s_k  = \text{min}\{ i\ | \ f(i,0) \leq k \}, $$
$$ \mu_i(k) = \text{min}\{j\ |\  f(i,j) \leq k \} \
\text{for} \ 0 \leq i \leq s_k - 1. $$
We call $\{\mu_i(k)\}$ the {\bf invariants} of $C$.
For $k >> 0$, $\mu_i(k) = \lm_i$, where the
$\{\lm_i \}$ are the invariants of a generic hyperplane
section, $\G = H \cap C$, of C.
The main result of this paper is to prove:

\begin{theorem}[The connectness of curve invariants] \label{theorem:connect}
Let $C$ be a reduced, irreducible, non-degenerate curve in $\pthree$,
then the invariants, $\{\mu_i(j)\}$, of $C$ are such that
$$ \mu_{i+1}(k)+2 \geq \mu_i(k) \geq \mu_{i+1}(k)+1 \
\text{for}\  0 \leq i < s_k-1, $$
and we say that $\{ \mu_i(k) \}$ is {\bf connected}.
Furthermore, if $s_k < s_0$, then $\mu_{s_k-1}(k) \leq 2$.
\end{theorem}

This result greatly restricts the kind of Borel-fixed monomial ideal which
can occur as the generic initial ideal of such curves
and thus gives us quite a lot of control over their Hilbert functions.

The paper is broken into parts as follows:
In Section 1, we will prove the Connectedness
Theorem. The proof incorporates the ideas used in the proof of
the Gruson, Peskine result ([GP], [EP]) and a more differential
approach due to Strano ([S]).
In Section 2, we will give some further restrictions
on the generic initial ideal of a reduced, irreducible,
non-degenerate curve in $\pthree$.
We will generalize a result of Strano on the
effect of certain generators of the generic initial ideal
of the curve on the syzygies of the
hyperplane section ideal to a larger group of ideals. We will also
generalize another result of Strano on complete intersections
to curves whose hyperplane sections have invariants of which the
first few are like those of a complete intersection. I believe
this last result will lead to some interesting results on the
gaps of Halphen ([E]).

{\it Acknowledgement.} This paper is taken from my dissertation and I
would like to take this opportunity to
thank my advisor Mark Green, for his insight and patience.

\section{Proof of the Connectedness Theorem}

\subsection{A pictorial description of monomial ideals}

We will first give an description, due to M. Green ([Gr]), of how
an ideal $I$ whose generators are monomials in
three variables may be described pictorially.
Hopefully, this desciption will make subsequent explanations clearer.
\newline
First, draw a triangle corresponding to the monomials
$x_1^ix_2^j$, where $i + j = n$, $0 \leq n \leq n_0$.
Let $f(i,j) = \text{min}\{k \ | \ x_1^ix_2^jx_3^k \in I \}$.
For each $i,j$,
if $f(i, j) = \infty $ (i.e. $x_1^ix_2^jx_3^k \notin I$ for all
$k \geq 0$) put a circle in the (i, j) position,
if $ 0 < f(i, j) < \infty$ put the number $f(i, j)$ in
the (i, j) position,
and if $f(i, j) = 0$ put an X in the (i, j) position.
\newline
Thus the triangle

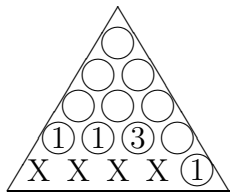
\vspace{5mm}
\begin{picture}(270,78)(0,0)

\put(150,75){\line(3,-5){42}}
\put(150,75){\line(-3,-5){42}}
\put(150,5){\line(1,0){42}}
\put(150,5){\line(-1,0){42}}

\put(150,61){\circle{12}}
\put(142.5,49){\circle{12}}
\put(157.5,49){\circle{12}}
\put(135,37){\circle{12}}
\put(150,37){\circle{12}}
\put(165,37){\circle{12}}
\put(127.5,25){\circle{12}}
\put(127.5,25){\makebox(0,0){1}}
\put(142.5,25){\circle{12}}
\put(142.5,25){\makebox(0,0){1}}
\put(157.5,25){\circle{12}}
\put(157.5,25){\makebox(0,0){3}}
\put(172.5,25){\circle{12}}
\put(120,13){\makebox(0,0){X}}
\put(135,13){\makebox(0,0){X}}
\put(150,13){\makebox(0,0){X}}
\put(165,13){\makebox(0,0){X}}
\put(180,13){\circle{12}}
\put(180,13){\makebox(0,0){1}}
\end{picture}

corresponds to the ideal
$(x_1^3x_3, x_1^2x_2x_3,  x_1x_2^2x_3^3, x_1^4, x_1^3x_2,
x_1^2x_2^2, x_1x_2^3, x_2^4x_3, x_2^5)$.
\newline
{\bf Note:} For (i, j) not in the picture, it is assumed
that $x_1^ix_2^j \in I$.
\newline
{\bf Remark:} It is the Borel-fixedness that ensures the step like look.
\newline
In this example the invariants are as follows:
$$ \begin{array}{llll}
\mu_0(0) = 5 & \mu_1(0) = 3 & \mu_2(0) = 2 & \mu_3(0) = 1\\
\mu_0(1) = 4 & \mu_1(1) = 3 & \mu_2(1) = 1 \\
\mu_0(2) = 4 & \mu_1(2) = 3 & \mu_2(2) = 1 \\
\mu_0(3) = 4 & \mu_1(3) = 2 & \mu_2(3) = 1
\end{array}$$
and $\mu_i(k) = \mu_i(3)$ for $k \geq 3$.

{\bf Note:} The invariants of this monomial ideal are connected.

To prove Theorem 1, we will first show in section 1.2, that
if there is an $i$ and $j$, with $0 \leq i < s_j-1$, such that
$\mu_{i+1}(j) + 2 < \mu_i(j)$, then
for generic hyperplanes $H, H'$, the ideal
$J = (I|_H:(H'\cap H)^j)$ can be ``split'' so that there exists
a homogeneous polynomial $X$ of degree $i+1$ such that
$\text{gin}(X)\cap\text{gin}(J) = \text{gin}(X \cap J)$.
We will then show in section 1.3, that if $I$ arises as the ideal of a
reduced, irreducible, non-degenerate  curve $C \in \pthree$,
such an $X$ cannot exist.

\subsection{Splitting a non-connected ideal}

{\bf Definition.}
An {\bf elementary move} $e_k$ for $ 1 \leq k \leq n$ is defined by:
$$ e_k(x^J) = x^{\hat{J}},$$
where
$\hat{J} = (j_1, \dots j_{k-1},j_k + 1, j_{k+1} - 1, j_{k+2}, \dots, j_{n+1})$
and where we adopt the convention that $x^J = 0 $ if $j_m < 0$ for
some $m$.

One can show that a monomial ideal $I$, is Borel-fixed if and only
if for all $x^J \in I$ and for every elementary move $e_k$,
$e_k(x^J) \in I$.

\begin{theorem}[syzygy configuration] \label{theorem:syzcon}
Let I be a Borel-fixed monomial ideal with generators
$x^{J_1}, \ \dots , x^{J_N}$. Then the first
syzygies of I is generated by
$$\{ x_i \otimes x^{J_j} - x^{L_{ij}} \otimes x^{J_{l_{ij}}} \ | \
1 \leq  j \leq N, \ 0 < i < max(J_j), \
\text{min}(L_{ij}) \geq \text{max}(J_{l_{ij}})\}, $$
where $\text{max}(J) = \text{max}\{ i \ |\  j_i > 0 \}$
and $\text{min}(J) = \text{min}\{i | j_i  > 0 \}$.
\end{theorem}

The proof of this is due to Green and may be found in [Gr].

\begin{theorem}\label{theorem:newgen}
Let $I$ be a homogeneous ideal generated in degree $\leq n$, with
generators $x^{J_1}, \ \dots , x^{J_N}$ of
$\gi$ in degree $\leq n$.
Then any generator $P$ of $\text{gin}(I)$ in degree $n+1$ is such that
$ P \prec  x_ix^{J_j}$, for some $J_j$ such that
$|J_j| = n$ and $ i < \text{max}(J_j)$.
\end{theorem}

{\bf Proof.}

We will first use the syzygies of $(x^{J_1}, \ \dots , x^{J_N})$
to obtain some new generators of $\text{gin}(I)_{n+1}$
which satisfy the conditions stated above.
\newline
$(x^{J_1}, \ \dots , x^{J_N})$ is a Borel fixed monomial ideal
so by Theorem~\ref{theorem:syzcon}, the
syzygies among the  $x^{J_j} $ are generated by syzygies  of the form
$$ \{\ x_k \otimes x^{J_j}- x^{L_{kj}} \otimes x^{J_{l_{kj}}} \ | \
1 \leq j \leq N, \ 1 \leq k < \text{max}(J_j)\}.$$
By Galligo's Theorem ([Ga]) we may assume,
after a possible change of basis, that $\gi = \text{in}(I)$.
Let $g_i \in I$ be such that $\text{in}(g_i) = x^{J_i}$ and assume the
$g_i$ are monic. Given a syzygy as above,
let $h_1 = x_kg_j-x^{L_{kj}}g_{l_{kj}}$,
then $\text{in}(h_1) \prec x_kx^{J_j}$.
\newline
Given $h_i$, if
$\text{in}(h_i)= x^{K_{i+1}}x^{J_{j_{i+1}}}$,
let $h_{i+1} = h_i - a_{i+1}x^{K_{i+1}}g_{j_{i+1}}$,
where $a_{i+1}$ is the leading coefficient of $h_i$, so that
$\text{in}(h_{i+1}) \prec \text{in}(h_i) \prec x_kx^{J_j}$.
\newline
This process must terminate, so for
$i$ sufficiently large we get either $h_i = 0$,
(in particular this will occur if deg$h_i \leq n$) or
$x^{J_j} $ does not divide $\text{in}(h_i)$ for all
$1 \leq j \leq N$, in which case
$\text{in}(h_i)$ is a new generator of $\gi$
satisfying the conditions stated in the theorem.
\newline
Now let $P = \text{in}(h) $ be a generator of $\text{gin}(I)_{n+1}$,
$h = \sum f_ig_i. $
\newline
Let $i_0$ be such that,
$\text{in}(f_{i_0}g_{i_0})$ is maximal, and
of the maximal $\text{in}(f_{i_0}g_{i_0})$,
$g_{i_0}$ has maximal degree, and among those $g_{i_0}$ of the same
maximal degree $\text{in}(g_{i_0})$ is minimal.
\newline
As $P \neq   \text{in}(f_{i_0}g_{i_0})$,
there exists $i_1$ such that
$\text{in}(f_{i_1}g_{i_1}) = \text{in}(f_{i_0}g_{i_0})$.
We have picked $i_0$  in such a way that either
deg$g_{i_0} > $ deg$g_{i_1}$, or
deg$g_{i_0} = $ deg$g_{i_1}$
and $\text{in}(g_{i_0}) \prec  \text{in}(g_{i_1})$.
\newline
We now want to show that $x_i | \text{in}(f_{i_0})$ for some
$ i < \text{max}(J_{i_0})$, where
$\text{in}(g_{i_0}) = x^{J_{i_0}}.$

{\bf Claim.}

Let $x^A, x^B$ be generating monomials of a Borel-fixed
monomial ideal, $I$, such that
$|A| > |B|$ and $x^Mx^A = x^Nx^B $ for some
monomials $x^M$ and $x^N$. Then there exists $x_i | x^M$
such that $i < \text{max}(A)$.

{\bf Proof of Claim.}

Let $A = (a_1, \dots ,a_s,0, \dots ,0)$, and
$B = (b_1, \dots ,b_s,b_{s+1}, \dots ,b_{n+1})$,
assume that $s = \text{max}(A)$.
Suppose $b_i \leq a_i$ for all $i \leq s$, then we may
apply elementary moves to $B$ to get $\hat B $
such that $x^{\hat{B}} \in I$ with
$\hat{b_i} \leq a_i$ for all $i$, and as $|\hat{B}| = |B| < |A| $,
this would imply $x^A$  cannot be a generator.
Therefore there exists $b_i > a_i$ for some $i \leq s$, and if
$i < s$ we are done. If $b_i \leq a_i$ for all $i < s$ but
$b_s > a_s$ we may again apply elementary moves to $B$ and get the
same conclusion. Therefore there exists $b_i > a_i $ for
$i < s$.

\vspace{1 mm}
If deg$g_{i_0} = $ deg$g_{i_1}$ and
$\text{in}(g_{i_0}) \prec \text{in}(g_{i_1})$, let
$\text{in}(g_{i_0}) = x^A$, $\text{in}(g_{i_1}) =  x^B$
where $A = (a_1, a_2, \dots a_{n+1})$ and $B = (b_1, b_2, \dots b_{n+1})$.
Then there exists $s$ such that
$a_k = b_k $ for all $k > s$ and $a_s > b_s$. As the degrees are
the same, there must exist $a_i < b_i$ for some $i < s$ and hence
$x_i | \text{in}(f_{i_0})$ for some
$ i <s  \leq  \text{max}(J_{i_0})$.
\newline
Thus, in either case we have $x_i | \text{in}(f_{i_0})$
with $ i < \text{max}(J_{i_0})$. Consider the syzygy
$$ x_i \otimes x^{J_{i_0}}- x^{L_{ii_0}} \otimes x^{J_{l_{ii_0}}}$$
and let $h^*$ be the element of $I$ obtained from this syzygy, as above.
\newline
$h^* =   x_ig_{i_0} - x^{L_{ii_0}}g_{l_{ii_0}} - \sum a_ix^{K_i}g_i$ where
$\text{in}(x^{K_i}g_i) \prec \text{in}(x_ig_{i_0})$,
and by Theorem~\ref{theorem:syzcon}
either $|J_{i_0}| > |J_{l_{ii_0}}|$
or $ \text{in}(g_{l_{ii_0}}) \succ \text{in}(g_{i_0})$.
Let $h_1 = h - e_{i_0} h^* $, where
$e_{i_0} = a_{i_0}\text{in}(f_{i_0})/x_i$ and $a_{i_0}$ is the leading
coefficient of $f_{i_0}$.
\newline
If deg$e_{i_0} = 0$,
$\text{in}(h_1) \preceq \text{max}\{\text{in}(h) ,\text{in}(h^*)\}$.
If in$(h) = $ in$(h^*)$, we are done as in$h^* \prec x_ix^{J_{i_0}}$.
Therefore we may assume
in$(h) \neq $ in$(h^*)$, and so either
in$(h_1) = $ in$(h^*)$ or in$(h_1) = $ in$(h)$. In the first case
we get $P = \text{in}(h) \prec \text{in}(h^*) \prec x_ix^{J_{i_0}}$,
and we are done.
\newline
If deg$e_{i_0} \geq 1$, $h^* = 0$ and
$\text{in}(h_1) = \text{in}(h) = P$.
\newline
Therefore we are left with the case in$(h_1) = P$, however in adding
the multiple of $h^*$, we only added terms which either have initial
term less than the term $\text{in}(f_{i_0}g_{i_0})$,
or if the term has the same initial term as the $g_{i_0}$ term,
corresponding to the $x^{L_{ii_0}}g_{l_{ii_0}}$ term, then either
$g_{l_{ii_0}}$ of lesser degree,
or  $\text{in}(g_{l_{ii_0}}) \succ \text{in}(g_{i_0}) $.
Therefore we may proceed by
induction and get the required result.\qed

Now suppose that $\gi$ is such that there exists
$i$ and $j$ such that
$$ \mu_{i+1}(j)+2 < \mu_i(j),
\ \text{where}\  0 \leq i <  s_j-1. $$
Consider the ideal $J = (I|_H:(H'\cap H)^j)$
for generic hyperplanes $H, H'$. Then
$\text{gin}(J) = (\text{gin}(I):x_3^j)$ and
$J$ has invariants $\nu_i(k) = \mu_i(k+j)$, and in particular
$$ \nu_{i+1}(0)+2 < \nu_i(0) . $$

For example, pictorially we could be in the following situation:

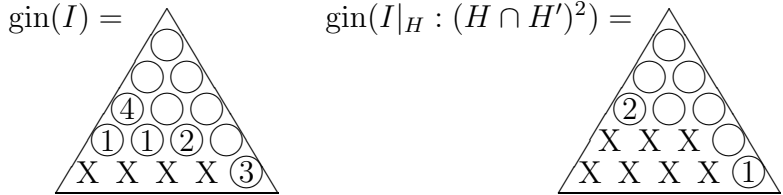
\vspace{.5cm}
\begin{picture}(300,78)(0,0)

\put(0,68){$\text{gin}(I) = $}
\put(60,75){\line(3,-5){42}}
\put(60,75){\line(-3,-5){42}}
\put(60,5){\line(1,0){42}}
\put(60,5){\line(-1,0){42}}

\put(60,61){\circle{12}}
\put(52.5,49){\circle{12}}
\put(67.5,49){\circle{12}}
\put(45,37){\circle{12}}
\put(45,37){\makebox(0,0){4}}
\put(60,37){\circle{12}}
\put(75,37){\circle{12}}
\put(37.5,25){\circle{12}}
\put(37.5,25){\makebox(0,0){1}}
\put(52.5,25){\circle{12}}
\put(52.5,25){\makebox(0,0){1}}
\put(67.5,25){\circle{12}}
\put(67.5,25){\makebox(0,0){2}}
\put(82.5,25){\circle{12}}
\put(30,13){\makebox(0,0){X}}
\put(45,13){\makebox(0,0){X}}
\put(60,13){\makebox(0,0){X}}
\put(75,13){\makebox(0,0){X}}
\put(90,13){\circle{12}}
\put(90,13){\makebox(0,0){3}}

\put(120,68){$\text{gin}(I|_H:(H\cap H')^2) = $}
\put(250,75){\line(3,-5){42}}
\put(250,75){\line(-3,-5){42}}
\put(250,5){\line(1,0){42}}
\put(250,5){\line(-1,0){42}}

\put(250,61){\circle{12}}
\put(242.5,49){\circle{12}}
\put(257.5,49){\circle{12}}
\put(235,37){\circle{12}}
\put(235,37){\makebox(0,0){2}}
\put(250,37){\circle{12}}
\put(265,37){\circle{12}}
\put(227.5,25){\makebox(0,0){X}}
\put(242.5,25){\makebox(0,0){X}}
\put(257.5,25){\makebox(0,0){X}}
\put(272.5,25){\circle{12}}
\put(220,13){\makebox(0,0){X}}
\put(235,13){\makebox(0,0){X}}
\put(250,13){\makebox(0,0){X}}
\put(265,13){\makebox(0,0){X}}
\put(280,13){\circle{12}}
\put(280,13){\makebox(0,0){1}}
\end{picture}

Let $K$ be the ideal generated by elements
of degree $\leq i + \nu_{i+1}(0) +2$ in $J$. We want to
show that there exists an ideal
$K \subseteq \hat{K} \subseteq J$ such that
$\text{gin}(\hat{K}) = (x_1^{i+1}) \cap \text{gin}(J)$.

\begin{corollary} \label{corollary:ngcor}
All elements of  $\text{gin}(K) $ are divisible by  $x_1^{i+1}$.
\end{corollary}

{\bf Proof.}

Let $x_1^ax_2^bx_3^c \in \text{gin}(K)_d$ for $d \leq i + \nu_{i+1}(0) + 2$.
If $a \leq i$ , then by Borel-fixedness
$x_1^ix_2^{\nu_{i+1}(0) + 2} \in \text{gin}(K) \subset \text{gin}(J)$, but
$x_1^ix_2^{\nu_i(0)} $ is a generator of
$\text{gin}(J) $ and so $\nu_{i+1}(0) + 2 \geq \nu_i(0)$.
But $ \nu_i(0) > \nu_{i+1}(0)+2 $.
\newline
Suppose by induction, all elements of  $\text{gin}(K)_d $ are
divisible by  $x_1^{i+1}$ for $d \geq i+\nu_{i+1}(0) + 2$.

{\bf Claim.}

If $d \geq  i+\nu_{i+1}(0) +2 $, then any generator of
$\text{gin}(K)_d$ has an $x_3$ term.

{\bf Proof of claim.}

If $x_1^ax_2^bx_3^c $ a gen\-er\-ator of $\text{gin}(K)_d$,
then $a \geq i+1$. Let $a = i+1+j$, we have
$x_1^{i+1+j}x_2^{\nu_{i+1+j}(0)} \in \text{gin}(J)$,
and $i+1+j + \nu_{i+1+j}(0) \leq i+1 + \nu_{i+1}(0) < d$,
hence $b < \nu_{i+1+j}(0)$ and $c > 0$.

\vspace{1 mm}
Let $P = x_1^ax_2^bx_3^c $ be a generator of $\text{gin}(K)_{d+1}$.
If  $a \leq i$, then by
Borel-fixedness $x_1^ix_2^{d+1-i} \in \text{gin}(K)_{d+1}$,
and as $x_1^{i+1} $ divides all elements of degree $\leq d$,
$x_1^ix_2^{d+1-i}$ is a generator of $\text{gin}(K)_{d+1}$,
and hence by Theorem~\ref{theorem:newgen}
$ x_1^ix_2^{d+1-i} \prec x_kx^{J}$ for
$x^{J}$ some generator of $\text{gin}(K)_d$
and $k < \text{max}(J)$. However, by the claim $x^{J}$ has an
$x_3$ term and thus
$x_1^ix_2^{d+1-i} \succ  x_kx^{J}$, which is a contradiction.
\qed

\begin{lemma}\label{lemma:lem}
Let $K$ be an ideal in $\C[x_1,x_2,x_3]$ such that
$\text{gin}(K) \subseteq (x_1^k)$, with $k \geq 1$ and
$k$ maximal.
Then, after a possible change of basis, there exists  a
homogeneous polynomial $X$ such that $in(X) = x_1^k$
and $K = XK'$.
\end{lemma}

We will need to use the following Proposition to prove
Lemma~\ref{lemma:lem} and later, to prove Theorem~\ref{theorem:connect}

\begin{proposition}\label{proposition:fam}

Let $\{ X_H \}$ be a family of homogeneous polynomials  in
the polynomial ring
$S = \C[x_1, \dots , x_{n+1}]$, (with $n \geq 2$)
parametrized by generic hyperplane sections $H= \sum t_ix_i$
in such a way that the coefficients of $X_H$ are differentiable
functions in the $t_i$.
Suppose that, if we differentiate with respect to
$H = \sum t_ix_i$ we get
$$x_j\frac{\p X_H}{\p t_i} - x_i\frac{\p X_H}{\p t_j}
\in (X_H) \ \text{mod}(H).$$
Then $X_H$ is a constant up to a multiple of $H$, for generic $H$.
\end{proposition}

{\bf Proof (Green).}

Let $Y = X_H|_H$, as
$(x_j\frac{\p X_H}{\p t_i} - x_i\frac{\p X_H}{\p t_j})|_H  \in (Y)$,
we may let
$(x_j\frac{\p X_H}{\p t_i} - x_i\frac{\p X_H}{\p t_j})|_H  = l_{ij}Y$.
$$x_k(x_j\frac{\p X_H}{\p t_i} - x_i\frac{\p X_H}{\p t_j}) -
x_j(x_k\frac{\p X_H}{\p t_i} - x_i\frac{\p X_H}{\p t_k}) +
x_i(x_k\frac{\p X_H}{\p t_j} - x_j\frac{\p X_H}{\p t_k}) = 0, $$
thus
$$ (x_kl_{ij} - x_jl_{ik} + x_il_{jk})Y = 0. $$
Assuming $Y \neq 0$ generically (otherwise we are done), we get
$$ x_kl_{ij} - x_jl_{ik} + x_il_{jk} = 0, $$
and hence
$$  l_{ij} = \al x_j + \be x_i, \  l_{ik} = \al x_k + \ga x_i, \
 l_{jk} = -\be x_k + \ga x_j,  $$
for some $\al$, $\be$ and $\ga$.
\newline
Therefore
$$(x_j(\frac{\p X_H}{\p t_i} - \al X_H)-
x_i(\frac{\p X_H}{\p t_j} + \be X_H))|_H = 0, $$
and we have
$$ \frac{\p X_H}{\p t_i} - \al X_H = x_iU \ \text{mod}(H),  \
\frac{\p X_H}{\p t_j} + \be X_H = x_jU \ \text{mod}(H)$$
and similarly
$$ \frac{\p X_H}{\p t_k} + \ga X_H = x_kU \ \text{mod}(H). $$
Letting $\al = \al_i$, $\be = -\al_j$, $\ga = -\al_k$ we get
$ \frac{\p X_H}{\p t_i} - x_iU = \al_iX_H \ \text{mod}(H) $.
We may vary $X_H$ continuously by any multiple of $H$ without changing the
hypothesis or claims of the Proposition, so letting
$X_H' = X_H - HU$ we get
$$
\frac{\p X_H'}{\p t_i} = \frac{\p X_H}{\p t_i} - x_iU - H\frac{\p U}{\p t_i}
 = \al_iX_H - H\frac{\p U}{\p t_i} \ \text{mod}(H) =
\al_iX_H' \ \text{mod}(H) $$
Therefore we may assume that $X_H$ is such that
$$ \frac{\p X_H}{\p t_i} = \al_iX_H + HU_i. $$
Differentiating again we get
$$\begin{array}{rl}
\frac{\p^2X_H}{\p t_j \p t_i} &=
\frac{\p \al_i}{\p t_j}X_H + \al_i\frac{\p X_H}{\p t_j} +
x_jU_i + H\frac{\p U_i}{\p t_j} \\
&= \frac{\p \al_i}{\p t_j}X_H + \al_i(\al_jX_H + HU_j)+
x_jU_i + H\frac{\p U_i}{\p t_j} \\
\end{array}$$
and
$$\begin{array}{rl}
\frac{\p^2X_H}{\p t_i \p t_j} &=
\frac{\p \al_j}{\p t_i}X_H + \al_j\frac{\p X_H}{\p t_i} +
x_iU_j + H\frac{\p U_j}{\p t_i} \\
&=  \frac{\p \al_j}{\p t_i}X_H + \al_j(\al_iX_H + HU_i) +
x_iU_j + H\frac{\p U_j}{\p t_i}.
\end{array}$$
Thus
$$ x_jU_i -x_iU_j = (\frac{\p \al_j}{\p t_i} - \frac{\p \al_i}{\p t_j})X_H
+ H(\al_jU_i - \al_iU_j + \frac{\p U_j}{\p t_i} - \frac{\p U_i}{\p t_j})$$
and so
$$( x_k(\frac{\p \al_j}{\p t_i} - \frac{\p \al_i}{\p t_j})
- x_j(\frac{\p \al_k}{\p t_i} - \frac{\p \al_i}{\p t_k})
+ x_i(\frac{\p \al_k}{\p t_j} - \frac{\p \al_j}{\p t_k}))X_H |_H = 0.$$
Therefore $\frac{\p \al_j}{\p t_i} = \frac{\p \al_i}{\p t_j}$ for all
$i, j$, and hence there exists $\al$ such that
$\al_i = \frac{\p \al}{\p t_i}$ for all $i$.
\newline
Thus we may assume
$ \frac{\p X_H}{\p t_i} = \frac{\p \al}{\p t_i}X_H + HU_i. $ Let
$X_H' = e^{-\al}X_H$, then
$$ \frac{\p X_H'}{\p t_i}
= e^{-\al}(\frac{\p X_H}{\p t_i}-\frac{\p \al}{\p t_i}X_H)
= 0 \ \text{mod}(H)$$
Therefore we may assume $ \frac{\p X_H}{\p t_i} = HU_i$.
\newline
Let $ \frac{\p X_H}{\p t_i} = H^kU_i$, with $k \geq 1$.
Let $I$ be such that $|I| = k$, then
$$ \frac{\p^{k+1}X_H}{\p t^I\p t_i} = k!x^IU_i \ \text{mod}(H).$$
If $I + v_i = I' + v_j$, then
$ \frac{\p^{k+1}X_H}{\p t^I\p t_i} = \frac{\p^{k+1}X_H}{\p t^{I'}\p t_j}$
and so $ (x^{I'}U_j - x^IU_i)|_H = 0$.
Now, $x^{I'} = x_ix^J$ and $x^I = x_jx^J$ and so
$(x_iU_j - x_jU_i)|_H = 0$, therefore
$U_i = x_iV + HV_i$
and  $ \frac{\p X_H}{\p t_i} = H^k(x_iV + HV_i)$.
Let
$X_H' = X_H - \frac{1}{k+1}H^{k+1}V$, then
$$ \begin{array}{rl}
\frac{\p X_H'}{\p t_i} &= \frac{\p X_H}{\p t_i}- H^kx_iV -
\frac{1}{k+1}H^{k+1}\frac{\p V}{\p t_i} \\
&= H^{k+1}V_i -\frac{1}{k+1}H^{k+1}\frac{\p V}{\p t_i} \\
&= H^{k+1}W_i.
\end{array}$$
Therefore we may assume $\frac{\p X_H}{\p t_i} = 0$ and
hence that $X_H$ is a constant up to a multiple of $H$.
\qed

{\bf Proof of Lemma~\ref{lemma:lem}}
\newline
$K \subseteq K^{sat}$, therefore
$\text{gin}(K) \subseteq \text{gin}(K^{sat})$ and
$\text{gin}(K^{sat})$ is generated by
$$\{ x_1^ax_2^b | x_1^ax_2^bx_3^c \in \text{gin}(K)
\ \text{for some} \  c \geq 0\}.$$
Therefore $\text{gin}(K^{sat}) \subseteq (x_1^k)$
so we may assume $K$ is saturated.
\newline Let $H$ be a generic hyperplane section,
then $\text{gin}((K|_H)^{sat}) = (x_1^k)$,
and hence  $(K|_H)^{sat} = (X_H)$.
\newline
Let $l$ be such that $(K|_H : \text{\bf m}^l) = (K|_H)^{sat} $.
Then for all monomials $x^I \in \text{\bf m}^l$,
there exists a polynomial $A_{H,I}$, such that
$x^IX_H + HA_{H,I} \in K$.
Differentiating with respect to $H = \sum t_ix_i$ we get
$$ x^I\frac{\p X_H}{\p t_i} + x_iA_{H,I} \in K|_H, $$
and thus
$$x_j\frac{\p X_H}{\p t_i} - x_i\frac{\p X_H}{\p t_j} \in
 (K|_H : \text{\bf m}^l) = (K|_H)^{sat}  = (X_H).$$
Thus considering the $X_H$ as elements  of $\C[x_1, x_2, x_3]$,
the $X_H$ satisfy the hypothesis of Proposition~\ref{proposition:fam}
and so $X = X_H|_H$
is independent of $H$.
\newline
Let $Y \in K$, Then $Y|_H \in K|_H \subseteq (K|_H)^{sat} = (X)$
and $Y = A_HX+HB_H \in K$. Differentiating with respect to $H$
we get
$0 =(\frac{\p A_H}{\p t_j}X + x_iB_H)|_H \in K|_H \subseteq (X)$.
Therefore $(x_iB_H)|_H \in (X)$ and thus $B_H|_H \in (X)$.
Let $B_H = XC_H + HB_H'$, then $Y = XA_H' + H^2B_H' \in (X)$.
\newline
Let $Y = A_HX+H^kB_H$, with $ k \geq 2$.
Differentiating with respect to $H$,
$k$ times, we get
$$0 =(\frac{\p^k A_H}{\p t^I}X + k!x^IB_H)|_H \in K|_H \subseteq (X).$$
Therefore $(x^IB_H)|_H \in (X)$ and thus $B_H|_H \in (X)$.
Letting $B_H = XC_H + HB_H'$, we get $Y = XA_H' + H^{k+1}B_H' \in (X)$.
Therefore we may assume $B_H = O$ and hence $Y = XA_H \in (X)$.
\qed

\noindent
Let $K$ be an ideal contained in $J$ maximal with
respect to the properties that
$K = XK'$ with $\text{in}(X) = x_1^{i+1}$,
and $(K)_d = (J)_d$ for $d \leq i+\nu_{i+1}(0) +2$.
$K = X \cap J$.

\begin{lemma}
$\text{gin}(X) \cap \text{gin}(J) \subseteq \text{gin}(X \cap J)$
\end{lemma}

{\bf Proof.}

By Galligo's Theorem ([Ga]), we may assume $\text{gin}(I) = \text{in}(I)$
for all ideals $I = (X), J$ or $K$.
\newline
Let $M = x_1^ax_2^bx_3^c \in \text{gin}(X) \cap \text{gin}(J)$,
then $a \geq i+1$, and we may write
$M = x_1^{a - \al}x_2^{b - \be}x_3^{c - \ga}( x_1^{\al}x_2^{\be}x_3^{\ga})$
where $A = x_1^{\al}x_2^{\be}x_3^{\ga}$ is a generator of
$\text{gin}(J)$.
If deg$A \leq i + \nu_{i+1}(0) + 2$, then
$A \in \text{gin}(X \cap J)$.
Suppose  deg$A > i + \nu_{i+1}(0) + 2$, and
$\al$ is maximal, if $\al < a$ and $\be $ or $\ga \geq 1$,
then either $x_1^{\al +1}x_2^{\be -1}x_3^{\ga}$ or
$x_1^{\al +1}x_2^{\be }x_3^{\ga - 1} \in \text{gin}(J)$. Then there
would exist $B = x_1^{\al'}x_2^{\be'}x_3^{\ga'} $, a generator of
$\text{gin}(J)$, dividing this element and $\al' = \al +1 $ otherwise
$A$ would not be a generator, and $B | M$ which contradicts the
maximality of $\al$. If
$\be = \ga = 0$, but then $x_1^{\al}$ is a generator of $J$
and so $\al \leq i + 1 + \nu_{i+1}$.
Therefore we may assume $\al = a \geq i+1$.
If $\ga \leq 1$ then deg$A \leq i+\nu_{i+1}(0) +2 $ therefore
we may assume $\ga \geq 2$.
\newline
Let's assume  $M = x_1^ax_2^bx_3^c \in \text{gin}(X) \cap \text{gin}(J)$,
is a generator of $\text{gin}(J)$ with deg$M > i+\nu_{i+1}(0)+2$,
$c \geq 2$, and suppose
$M \not\in \text{gin}(X \cap J) = \text{gin}(K)$.
Pick $M$ satisfying this of minimal degree $m$, and among those
of minimal degree, let $M$ be maximal. Pick $f \in J$
with $\text{in}(f) = M$.
Let  $L$ be the ideal generated
by $K$ and $f$. Then $\text{gin}(L)_d = \text{gin}(K)_d$ for $d < m$
and $\text{gin}(L)_m = \text{gin}(K)_m + M$.
Then the generators of $\text{gin}(L)$ in degree $\leq m$ form
a Borel-fixed monomial ideal and every generator of degree $m$
has an $x_3$ term.
Thus as in the proof of
Theorem~\ref{theorem:newgen} and Corollary~\ref{corollary:ngcor}
all elements of $\text{gin}(L)$ must be
divisible by $x_1^{i+1}$. This however contradicts the maximality
of $K$, therefore
$\text{gin}(X) \cap \text{gin}(J) \subset \text{gin}(X \cap J)$.
\qed

\subsection{Proof of the Connectedness Theorem}

\noindent
Suppose $\mu_{i+1}(j) + 2 < \mu_{i}(j)$ for $i < s_j -1$.
Consider $J = (I|_H:(H'\cap H)^j)$,
for generic hyperplanes $H, H'$ and let
$K$ be the ideal contained in $J$
such that $(K)_d = (J)_d$ for $d \leq i+2+\mu_{i+1}(j)$ and
$K = XK'$ with $\text{in}(X) = x_1^{i+1}$.

\begin{proposition}
$X$ is constant up to a multiple of $H$.
\end{proposition}

{\bf Proof.}

\noindent
Pick $ p = p(H,H') \in K'$ such that
deg$(pX) = m \leq i+1+\mu_{i+1}(j)$ and $(p,X) = 1$.
(As $x_1^{i+1}x_2^{\mu_{i+1}(j)} \in
\text{gin}(J)_{i+1+\mu_{i+1}(j)} = \text{gin}(K)_{i+1+\mu_{i+1}(j)}$,
there exists $p \in K'$ such that in$(p) = x_2^{\mu_{i+1}(j)}$,
and so $p$ and $X$ cannot have a common factor.)
\newline
Then
$$ (H')^jpX + HA \in I \ \text{for some} \ A = A(H,H').$$
Keeping $H'$ fixed and letting $H = \sum t_ix_i$ vary and
differentiating with respect to  the $t_i$ we get
$$ (H')^j(\frac{\p p}{\p t_i}X + p\frac{\p X}{\p t_i}) + x_iA \in I_{H}$$
and so
$$ x_k(\frac{\p p}{\p t_i}X + p\frac{\p X}{\p t_i}) -
x_i(\frac{\p p}{\p t_k}X + p\frac{\p X}{\p t_k})
\in (I_H : (H')^j)_{m+1} = J_{m+1}.$$
$m+1  \leq i+2+\mu_{i+1}(j)$, hence
$(J)_{m+1} = (K)_{m+1} \subseteq (X)$.
Therefore $p(x_j\frac{\p X}{\p t_i} - x_i\frac{\p X}{\p t_j}) \in (X)$,
$(p, X) = 1$ and so
$x_j\frac{\p X}{\p t_i} - x_i\frac{\p X}{\p t_j} \in (X)$.
Therefore, for generic $H$, $X_H = X$, satisfy the hypothesis
of Proposition~\ref{proposition:fam}
and hence we may assume $X$ is constant up to
a multiple of $H$.
\qed

{\bf Proof of Theorem~\ref{theorem:connect}}

Keeping $H'$ fixed and letting $Y = (H')^jX$, we have for
all hyperplanes $H$ and all $p_H =p(H,H') \in K'_{H,H'}$ there exists
$A_H = A(H,H')$ such that
$$ p_HY + HA_H \in I. $$
Let $\G_H = H \cap C$, then
$\G_H \subset V(p_HY) = V(Y) \cup V(p_H)$.
If for a generic $H$, there exists $q_H $ a point in $\G_H$ such that
$q_H \in V(Y)$, then
$S = \{q_H \ | \ q_H \in V(Y) \}$ will be a 1-dimensional
space and $S \subset C$, $C$ is reduced and irreducible
therefore $S$ is dense in  $C$ and hence $ C \subset V(Y)$.
However $V(Y) = V((H')^j) \cup V(X)$ and as
C is nondegenerate, this would imply $ C \subset V(X)$.
However, $i < s_j - 1 \leq s_0 - 1$ and so $i+1 < s_0$.
But $s_0$ is the smallest degree of elements of $I$,
therefore $C \not\subset V(X)$.
Therefore for generic $H, H'$,  $ \G_H \subset V(K')$.
However the invariants of $(K')^{sat}$ are
$\lm_{i+1} > \dots > \lm_{s-1}$ and as $i+1 > 0$,
$ \G_H \not\subset V(K')$.
This concludes the proof of the main part of Theorem~\ref{theorem:connect}.

\vspace{2mm}
Now suppose $s_k < s_0$, and $\mu_{s_k-1}(k) \geq 3$.
For generic hyperplane sections $H$ and $H'$,
consider $(I|_H : (H'\cap H)^k)$,
we have $(I|_H : (H'\cap H)^k)_{s_k} = (X_{H,H'})_{s_k}$, and
$(I|_H : (H'\cap H)^k)_{s_k+1} = (X_{H,H'})_{s_k+1}$, and so
there exists
a homogeneous polynomial $A_{H,H'}$ such that
$$ (H')^kX_{H,H'} + HA_{H,H'} \in I.$$
Keeping $H'$ fixed and
differentiating with respect to $H = \sum t_ix_i$, we get
$$ (H')^k(\frac{\p X_{H,H'}}{\p t_i}) + x_iA_{H,H'} \in I_{H}$$
and so
$$ (x_j\frac{\p X_{H,H'}}{\p t_i} -
x_i\frac{\p X_{H,H'}}{\p t_j})
\in (I|_H : (H'\cap H)^k)_{s_k+1} = (X_{H,H'})_{s_k+1}.$$
Thus by Proposition~\ref{proposition:fam},
$X_{H,H'}$ is constant with respect to
$H$.
\newline
Let $X_{H,H'} = X_{H'}$, Fix $H'$ and let
$\G_{H}$ be a generic hyperplane section of $C$.
As $ (H')^kX_{H'} + HA_{H,H'} \in I$,
$\G_{H} \subset V(H') \cup V(X_{H'})$.
But as the points of $\G_{H}$ are in general position,
there must exist at least one point of $\G_{H} \in V(X_{H'})$.
But varying $H$ as above would imply that
$C \subset V(X_{H'})$, and
hence that $s_k \geq s_0$, which is a contradiction.
This completes the proof of Theorem~\ref{theorem:connect}. \qed

\section{Further Results on the Generic Initial
Ideal of a Curve}

\subsection{Generalized Strano}

\noindent
This result generalizes a result of Strano ([S]):

\noindent
{\bf Definition.}
If $x_1^ix_2^jx_3^{f(i,j)}$ is a generator of $\text{gin}(I_C)$ with
$f(i,j) > 0$, then $x_1^ix_2^jx_3^k$ is a {\bf sporadic zero} for
all $0 \leq k < f(i,j)$.

\begin{theorem}[Strano] \label{theorem:strano}
If C is a reduced irreducible curve and has a sporadic zero
in degree m, then $I_{\G}$ has a syzygy in degree $\leq m+2$.
\end{theorem}

\begin{theorem}[Generalized Strano] \label{theorem:genstr}
Let C be a reduced irreducible curve with a sporadic zero
$x_1^ix_2^jx_3^{k-a}$ of degree m, such that
$x_1^ix_2^jx_3^k$ is a generator of  $\text{gin}(I_C)$.
Then, for generic hyperplanes $H$ and $H'$,
$J = (I_C|_H : (H \cap H')^a)$ has a syzygy in degree $\leq m+2$.
\end{theorem}

\noindent{\bf Proof.}

\noindent
$x_1^ix_2^jx_3^{k-a} \in \text{gin}(J)_m$, therefore
there exists $F \in (I_C|_H : (H \cap H')^a)_m = (J)_m$
varying  differentiably with $H$ and $H'$, and hence
for generic $H$ and $H'$ there exists $A$, depending on
$H$ and $H'$ such that
$$ (H')^aF + HA \in I_C.$$
If we were to keep $H'$ fixed, the coefficients of
$F$ and $A$ would be homogeneous polynomials in the
coefficients $t_i$ of $H = \sum t_ix_i$. We will choose
$F$, which is a bihomogeneous polynomial in $t_i$ and
$x_i$, such that the degree of $F$ with respect to
$t_i$ is minimal. We will also assume that
$\text{gin}(J) = \text{in}(J)$ and that
$\text{in}(F) = x_1^ix_2^jx_3^{k-a}$.
\newline
Differentiating with respect to $H$, keeping $H'$ fixed we
get
$$ (H')^a\frac{\p F}{\p t_j} + x_jA \in I_C|_H$$
and so
$$x_j\frac{\p F}{\p t_i} - x_i\frac{\p F}{\p t_j} \in (J)_{m+1}.$$
Hence
$$x_k(x_j\frac{\p F}{\p t_i} - x_i\frac{\p F}{\p t_j}) -
x_j(x_k\frac{\p F}{\p t_i} - x_i\frac{\p F}{\p t_k}) +
x_i(x_k\frac{\p F}{\p t_j} - x_j\frac{\p F}{\p t_k}) = 0 $$
is a syzygy of $J$ in degree $m+2$.

\noindent
Suppose $J$ does not have a syzygy in degree $\leq m+2$, then
$$x_j\frac{\p F}{\p t_i} - x_i\frac{\p F}{\p t_j}
= x_jU_i - x_iU_j \ \text{where} \ U_i \in (J)_m.$$
Rewriting, we get
$$x_j(\frac{\p F}{\p t_i}-U_i) - x_i(\frac{\p F}{\p t_j}-U_j) = 0$$
and so
$$\frac{\p F}{\p t_i} = U_i  + x_iR.$$
Letting $F' = F- HR$ we get
$$
\frac{\p F'}{\p t_i} = \frac{\p F}{\p t_i} - x_iR - H\frac{\p R}{\p t_i}
 = U_i \ \text{mod}(H) \in (J)_m.$$
As we have assumed the degree of $F$ is minimal with respect to
$t_i$ we get that $F$ is constant up to a multiple of $H$.
Hence, by an argument similar to that of Theorem~\ref{theorem:connect},
$F \in I_C$. This, however, is a contradiction.
\qed

\noindent
{\bf Example.}
\newline
The following diagram can not correspond to a generic initial
ideal of a curve, even though it is connected.

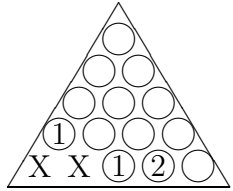
\vspace{5mm}

\begin{picture}(370,75)(0,7)

\put(100,75){\line(3,-5){42}}
\put(100,75){\line(-3,-5){42}}
\put(100,5){\line(1,0){42}}
\put(100,5){\line(-1,0){42}}

\put(100,61){\circle{12}}
\put(92.5,49){\circle{12}}
\put(107.5,49){\circle{12}}
\put(85,37){\circle{12}}
\put(100,37){\circle{12}}
\put(115,37){\circle{12}}
\put(77.5,25){\circle{12}}
\put(77.5,25){\makebox(0,0){1}}
\put(92.5,25){\circle{12}}
\put(107.5,25){\circle{12}}
\put(122.5,25){\circle{12}}
\put(70,13){\makebox(0,0){X}}
\put(85,13){\makebox(0,0){X}}
\put(100,13){\circle{12}}
\put(100,13){\makebox(0,0){1}}
\put(115,13){\circle{12}}
\put(115,13){\makebox(0,0){2}}
\put(130,13){\circle{12}}

\end{picture}

\noindent
We have a sporadic zero in degree $3$, and so by the Theorem,
$J = (I_C|_H : (H)^1)$ has a syzygy in degree $\leq 5$.
The diagram of $\text{gin}(J)$ is as follows:

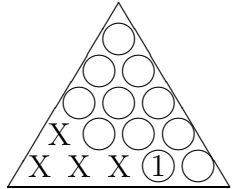
\begin{picture}(370,75)(0,7)

\put(100,75){\line(3,-5){42}}
\put(100,75){\line(-3,-5){42}}
\put(100,5){\line(1,0){42}}
\put(100,5){\line(-1,0){42}}

\put(100,61){\circle{12}}
\put(92.5,49){\circle{12}}
\put(107.5,49){\circle{12}}
\put(85,37){\circle{12}}
\put(100,37){\circle{12}}
\put(115,37){\circle{12}}
\put(77.5,25){\makebox(0,0){X}}
\put(92.5,25){\circle{12}}
\put(107.5,25){\circle{12}}
\put(122.5,25){\circle{12}}
\put(70,13){\makebox(0,0){X}}
\put(85,13){\makebox(0,0){X}}
\put(100,13){\makebox(0,0){X}}
\put(115,13){\circle{12}}
\put(115,13){\makebox(0,0){1}}
\put(130,13){\circle{12}}

\end{picture}

\noindent
$J$ has only two generators in degree $\leq 4$, corresponding to
$x_1^3$ and $x_1^2x_2^2$ and so if there were a syzygy in degree
$\leq 5$ this would imply that we may ``split'' the ideal as in the
proof of Theorem~\ref{theorem:connect}. But this would give a contradiction.

\vspace{.5cm}
\noindent
More generally, if $s_k = \text{min}\{i | f(i,0) \leq k\}$ and
for $k >> 0 $ $s_0 > s_k $ and $\mu_{s_k-1}(k) = \lambda_{s-1} = 2$,
then $f(s-2,3) \leq f(s,0)$.

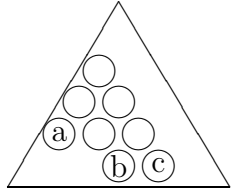
\begin{picture}(370,75)(0,7)

\put(100,75){\line(3,-5){42}}
\put(100,75){\line(-3,-5){42}}
\put(100,5){\line(1,0){42}}
\put(100,5){\line(-1,0){42}}

\put(92.5,49){\circle{12}}
\put(85,37){\circle{12}}
\put(100,37){\circle{12}}
\put(77.5,25){\circle{12}}
\put(77.5,25){\makebox(0,0){a}}
\put(92.5,25){\circle{12}}
\put(107.5,25){\circle{12}}
\put(100,13){\makebox(0,0){b}}
\put(100,13){\circle{12}}
\put(115,13){\circle{12}}
\put(115,13){\makebox(0,0){c}}

\end{picture}

By connectedness $b \leq a$. If $c > a$, then if we let
$J = (I_C|_H : (H)^a)$, then $J$ has a syzygy in degree $\leq s+2$.
But again this would imply that we could ``split'' the ideal.

\subsection{Complete Intersections and Almost Complete Intersections}

The result in this section is inspired by the work of Ellia ([E]) and
again generalizes a result of Strano ([S]):

\begin{theorem}[Strano]
If $C$ is a reduced irreducible curve whose generic hyperplane
section has the Hilbert function of a complete intersection
of type $(m,n)$, where $m,n > 2$, then $C$ is a complete
intersection of type $(m,n)$.
\end{theorem}

This result follows from connectedness and Theorem~\ref{theorem:strano}.

\noindent
{\bf Note.}
If $\G$ is a set of $d$ points in general position with invariants
$\lm_0 > \dots > \lm_{k-1} > 0$ such that $\lm_i = \lm_0 - 2i$ for all $i$.
Then $\G$ is a complete intersection of type $(k, d/k)$. (See [Gr].)

\begin{proposition}
Let C be a reduced, irreducible, non-degenerate curve in $\pthree$, let
$\Gamma = C \cap H$ be a generic hyperplane section with invariants
$\{\lambda_i\}_{i=0}^{s-1}$. If $\lambda_{s-i} = \lambda_{s-1} +2(i-1)$
for $1 \leq i \leq k$, where $k \geq 3$, then $f(i,j) > 0 $ only
if $i < s-k$.
\end{proposition}

\noindent
{\bf Proof.}
\newline
Let $J = (I|_H : (H \cap H')^j)$ for $j >> 0$, so that
$\text{gin}(J) = \text{gin}(I_{\Gamma})$. Let $f$ correspond
to $x_1^s$ and let $g$ correspond to $x_1^{s-1}x_2^{\lambda_{s-1}}$,
where $f$ and $g$ are in $J$. If $f$ and $g$ have a syzygy in
degree $d \leq \lambda_{s-k} + (s-k)$, then generators of
$\text{gin}(J)$ in degree $d$ correspond to generators of $J$ and thus
we may ``split'' the ideal $J$ as in the proof of
Theorem~\ref{theorem:connect}.
Therefore $f$ and $g$ have no syzygy in degree $\leq \lambda_{s-k} + (s-k)$.
By Theorem~\ref{theorem:strano}  or Theorem~\ref{theorem:genstr},
this means that there can be no
sporadic zeroes in degree $\leq \lambda_{s-k} + (s-k) -2$.
\newline
If there is a sporadic zero in degree $\lambda_{s-k} + (s-k) -1 = $
$\lambda_{s-(k-1)} +(s-(k-1)$, then $\mu_{s-(k-1)}(0) > \lambda_{s-(k-1)}$
and $\mu_{s-(k-2)}(0) = \lambda_{s-(k-2)} = \lambda_{s-(k-1)} - 2$,
which contradicts the connectedness of the $\{\mu_i(0)\}$.
Similarly if $f(s-k, \lambda_{s-k}) > 0$ then
$\mu_{s-k}(0) > \lambda_{s-k} = \lambda_{s-(k-1)} + 2 =
\mu_{s-(k-1)}(0) + 2$ which again contradicts connectedness.
\qed

\noindent
Thus for the following configuration of a hyperplane section, we can only
possibly get a sporadic zero in the (0,6) position.

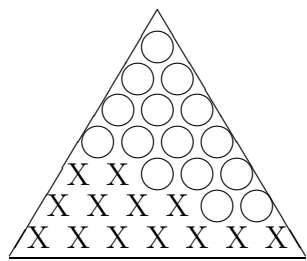
\vspace{5mm}

\begin{picture}(370,99)(0,0)

\put(140,99){\line(3,-5){56}}
\put(140,99){\line(-3,-5){56}}
\put(140,5){\line(1,0){56}}
\put(140,5){\line(-1,0){56}}

\put(140,85){\circle{12}}
\put(132.5,73){\circle{12}}
\put(147.5,73){\circle{12}}
\put(125,61){\circle{12}}
\put(140,61){\circle{12}}
\put(155,61){\circle{12}}
\put(117.5,49){\circle{12}}
\put(132.5,49){\circle{12}}
\put(147.5,49){\circle{12}}
\put(162.5,49){\circle{12}}
\put(110,37){\makebox(0,0){X}}
\put(125,37){\makebox(0,0){X}}
\put(140,37){\circle{12}}
\put(155,37){\circle{12}}
\put(170,37){\circle{12}}
\put(102.5,25){\makebox(0,0){X}}
\put(117.5,25){\makebox(0,0){X}}
\put(132.5,25){\makebox(0,0){X}}
\put(147.5,25){\makebox(0,0){X}}
\put(162.5,25){\circle{12}}
\put(177.5,25){\circle{12}}
\put(95,13){\makebox(0,0){X}}
\put(110,13){\makebox(0,0){X}}
\put(125,13){\makebox(0,0){X}}
\put(140,13){\makebox(0,0){X}}
\put(155,13){\makebox(0,0){X}}
\put(170,13){\makebox(0,0){X}}
\put(185,13){\makebox(0,0){X}}

\end{picture}

\section*{\bf References}

{\bf [B]} D. Bayer, {\it The division algorithm and
the Hilbert scheme}, Ph.D. Thesis, Harvard University,
Department of Mathematics, June 1982. Order Number 82-22588,
University Microfilms International, 300N Zeeb Rd., Ann Arbor,
MI 48106.

{\bf [BM]} D. Bayer, D. Mumford {\it What can be computed in
Algebraic Geometry}, Computational Algebraic Geometry and
Commutative Algebra, Symposia Mathematica Volume XXXIV, Cambridge
University Press. (1991) pp1-48.

{\bf [BS]} D. Bayer, M. Stillman {\it A criterion for
detecting m-regularity}, Invent. Math. 87 (1987) pp1-11.

{\bf [E]} P. Ellia {\it Sur les lacunes d'Halphen}, Algebraic curves and
projective geometry (Trento, 1988), 43-65, Lecture Notes in Math.,
1389, Springer, Berlin-New York 1989.

{\bf [EP]} P. Ellia, C. Peskine {\it Groupes de points de $\ptwo$:
caract\`ere et position uniforme}, Algebriac Geometry
(L'Aquila 1988), Lecture Notes in
Mathematics, 1417, Springer, Berlin 1990 pp111-116.

{\bf [Ga]} A. Galligo {\it A propos du theoreme de preparation de
Weierstrass}, Fonctions de Plusieurs Variables Complexes, Lecture Notes
in Math., Vol 1974 pp543-579.

{\bf [Gr]} M. Green {\it Lecture notes on generic initial ideals},
Unpublished.

{\bf [GP]} L. Gruson, C. Peskine {\it Genre des courbes de
l'espace projectifs}, Lecture Notes in Mathematics, 687,
Springer (1978) pp31-59.

{\bf [S]} R. Strano {\it Sulle Sezione Iperpiane Delle Curve},
Rend. Sem. Mat. Fis. Milano 57 (1987) pp125-134.

\vspace{.5cm}
Michele Cook

Department of Mathematics

UCLA

Los Angeles, CA. 90024

e-mail shelly\verb+@+math.ucla.edu

\end{document}